\documentclass[12pt]{article}
\usepackage{epsfig}
\usepackage{exscale}
\usepackage{fullpage}
\usepackage{url}

\usepackage{amsfonts}
\newcommand{\CC}{{\mathbb C}}
\newcommand{\HH}{{\mathbb H}}
\newcommand{\KK}{{\mathbb K}}
\newcommand{\OO}{{\mathbb O}}
\newcommand{\RR}{{\mathbb R}}

\newcommand{\JJ}{{\bf H}}
\newcommand{\mm}{{\bf M}}

\newcommand{\tsl}{TSL}
\newcommand{\ssl}{SL}
\newcommand{\ssla}{SL_1}
\newcommand{\sslb}{SL_2}

\renewcommand\AA{{\cal A}}
\newcommand{\II}{{\cal I}}
\newcommand{\MM}{{\cal M}}
\newcommand{\PP}{{\cal P}}
\newcommand{\QQ}{{\cal Q}}
\newcommand{\TT}{{\cal T}}
\newcommand{\XX}{{\cal X}}
\newcommand{\YY}{{\cal Y}}
\newcommand{\ZZ}{{\cal Z}}

\newcommand{\nxn}{$n \times n$}
\renewcommand{\bar}{\overline}
\newcommand{\tr}{{\rm tr}}
\renewcommand{\Re}{{\rm Re}}

\newcommand{\ibar}{\bar\imath}
\newcommand{\jbar}{\bar\jmath}
\newcommand{\kbar}{\bar k}

\newcommand{\spinor}{Cayley spinor}
\newcommand{\Spinor}{Cayley spinor}

\let\tilde=\widetilde

\begin{document}


\title{
{\bfseries\boldmath Octonionic Cayley Spinors and $E_6$}
}

\author{
	Tevian Dray \\[-2.5pt]
	\normalsize
	\textit{Department of Mathematics, Oregon State University,
		Corvallis, OR  97331} \\[-2.5pt]
	\normalsize
	{\tt tevian{\rm @}math.oregonstate.edu} \\
		\and
	Corinne A. Manogue \\[-2.5pt]
	\normalsize
	\textit{Department of Physics, Oregon State University,
		Corvallis, OR  97331} \\[-2.5pt]
	\normalsize
	{\tt corinne{\rm @}physics.oregonstate.edu}
}

\date{\normalsize February 14, 2010}

\maketitle

\begin{abstract}
Attempts to extend our previous work using the octonions to describe
fundamental particles lead naturally to the consideration of a particular
real, noncompact form of the exceptional Lie group~$E_6$, and of its
subgroups.  We are therefore led to a description of~$E_6$ in terms of
$3\times3$ octonionic matrices, generalizing previous results in the
$2\times2$ case.  Our treatment naturally includes a description of several
important subgroups of~$E_6$, notably $G_2$, $F_4$, and (the double cover of)
$SO(9,1)$.  An interpretation of the actions of these groups on the squares of
3-component \textit{\spinor s} is suggested.
\end{abstract}


\section{Introduction}

In previous work~\cite{Dim,Spin}, we used a formalism involving $2\times2$
octonionic matrices to describe the Lorentz group in 10 spacetime dimensions,
and then applied this formalism to the Dirac equation.  We developed a
mechanism for reducing 10 dimensions to 4 without compactification, thus
reducing the 10-dimensional massless Dirac equation to a unified treatment of
massive and massless fermions in 4~dimensions.  This description involves both
vectors (momentum) and spinors (solutions of the Dirac equation), which we
here combine into a single, 3-component object.  This leads to a
representation of the Dirac equation in terms of $3\times3$ octonionic
matrices, revealing a deep connection with the exceptional Lie group $E_6$.

\section{The Lorentz Group}

In earlier work~\cite{Lorentz}, we gave an explicit octonionic representation
of the finite Lorentz transformations in 10 spacetime dimensions, which we now
summarize in somewhat different language.

Matrix groups are usually defined over the complex numbers $\CC$, such as the
Lie group $SL(n;\CC)$, consisting of the \nxn\ complex matrices of determinant
$1$, or its subgroup $SU(n;\CC)$, the unitary (complex) matrices with
determinant $1$.  It is well-known that $SL(2,\CC)$ is the the double cover of
the Lorentz group $SO(3,1)$ in 4 spacetime dimensions, $\RR^{3+1}$.  One way
to see this is to represent elements of $\RR^{3+1}$ as $2\times2$ complex
Hermitian matrices $X\in \JJ_2(\CC)$, noting that $\det X$ is just the
Lorentzian norm.  Elements $M\in SL(2;\CC)$ act on $X\in \JJ_2(\CC)$ via
linear transformations of the form
\begin{equation}
T_M(X) = MXM^\dagger
\label{2action}
\end{equation}
and such transformations preserve the determinant.  The set of transformations
of the form~(\ref{2action}) with $M\in SL(2;\CC)$ is a group under
composition, and is therefore isomorphic to, and can be identified with,
$SO(3,1)$.
However, the map
\begin{eqnarray}
SL(2;\CC) &\longrightarrow& SO(3,1)\\\nonumber
M &\longmapsto& T_M
\end{eqnarray}
which takes $M$ to the linear transformation defined by~(\ref{2action}), is
not one-to-one; in fact, this map is easily seen to be a two-to-one
homomorphism with kernel $\{\pm I\}$.  We call such a homomorphism a
\textit{double cover}.  Restricting $M$ to the subgroup $SU(2;\CC)\subset
SL(2;\CC)$ similarly leads to the well-known double cover
\begin{eqnarray}
SU(2;\CC) &\longrightarrow& SO(3)
\end{eqnarray}
of the rotation group in three dimensions.  It is straightforward to restrict
the maps above to the reals, obtaining the double covers
\begin{eqnarray}
	SL(2;\RR) &\longrightarrow& SO(2,1) \\
	SU(2;\RR) &\longrightarrow& SO(2)
\end{eqnarray}

Since determinants of non-Hermitian matrices over the division algebras $\HH$
and $\OO$ are not well-defined, we seek alternative characterizations of these
complex matrix groups which do not involve such determinants.  The key idea is
that the determinant of ($2\times2$ and $3\times3$) Hermitian matrices over
any division algebra $\KK=\RR,\CC,\HH,\OO$ is well-defined, and therefore so
is the notion of determinant-preserving transformations.
We therefore define
\begin{equation}
\tsl(2;\HH) := \left\{
	T_M : \det\bigl(T_M(X)\bigr)=\det X \quad\forall X\in\JJ_2(\HH)\right\}
\label{tsldef}
\end{equation}
to be the set of determinant-preserving transformations in the quaternionic
case, where $M$ is now a quaternionic $2\times2$ matrix.  It is
straightforward to verify that $\tsl(2;\HH)$ is a group under composition, and
that
\begin{equation}
\tsl(2;\HH) \cong SO(5,1)
\end{equation}
under which we identify quaternionic linear transformations of the
form~(\ref{2action}) with the corresponding Lorentz transformations in
$\RR^{5+1}$.

We also have the \textit{spinor action} of $2\times2$  quaternionic matrices
$M$ on 2-component column vectors, namely
\begin{equation}
S_M(v) = Mv
\label{spinor}
\end{equation}
with $v\in\HH^2$.  We now \textit{define} $SL(2;\HH)$ to be the spinor
transformations $S_M$ such that the corresponding (vector) transformation
$T_M$ is determinant-preserving, that is,
\begin{equation}
\ssl(2;\HH) := \left\{ S_M : T_M\in\tsl(2;\HH) \right\}
\end{equation}
and it is straightforward to verify that this set of linear transformations
is a group under composition.  Furthermore, the map
\begin{eqnarray}
\ssl(2;\HH) &\longrightarrow& \tsl(2;\HH) \\\nonumber
S_M &\longmapsto& T_M
\end{eqnarray}
is again easily seen to be a two-to-one homomorphism, this time with kernel
$\{S_{\pm I}\}$, leading to the double cover
\begin{eqnarray}
	SL(2;\HH) &\longrightarrow& SO(5,1)
\end{eqnarray}
Requiring in addition that $\tr(MXM^\dagger) = \tr\ X$ for all
$X\in\JJ_2(\HH)$, and repeating the above construction, leads to the subgroup
$SU(2;\HH)\subset SL(2;\HH)$ and the double cover
\begin{eqnarray}
	SU(2;\HH) &\longrightarrow& SO(5)
\end{eqnarray}

Generalizing these groups to $\OO$ must be done with some care due to the lack
of associativity; for this reason, most authors discuss the corresponding Lie
algebras instead.  However, since composition of transformations of the
forms~(\ref{2action}) or~(\ref{spinor}) is associative, the above construction
can indeed be generalized~\cite{Lorentz}, provided care is taken
that~(\ref{2action}) itself is well-defined, that is, provided we require $M$
to satisfy
\begin{equation}
M(XM^\dagger) = (MX)M^\dagger
\label{welldef}
\end{equation}
for all $X\in\JJ_2(\OO)$.  In order to be able to later combine spinor
transformations $S_M$ with vector transformations $T_M$, we also require our
transformations to be \textit{compatible}~\cite{Lorentz,Mobius} with the
mapping from spinors to vectors given by $v\mapsto vv^\dagger$.  Explicitly,
we require
\begin{equation}
S_M(v)\bigl(S_M(v)\bigr)^\dagger = T_M(vv^\dagger)
\end{equation}
or in other words
\begin{equation}
(Mv)(v^\dagger M^\dagger) = M(vv^\dagger)M^\dagger
\label{comp}
\end{equation}
for all $v\in\OO^2$.  Conditions~(\ref{welldef}) and~(\ref{comp}) turn out to
be equivalent to the assumption that $M$ is complex
\footnote{A \textit{complex matrix} is one whose elements lie in a complex
subalgebra of the division algebra in question, in this case $\OO$.  Each such
matrix has a well-defined determinant.  It is important to note that there is
no requirement that the elements of two such matrices lie in the \textit{same}
complex subalgebra.}
and that
\begin{equation}
\det M \in \RR
\end{equation}
We therefore let $\mm_2(\OO)$ denote the set of \textit{complex} $2\times2$
octonionic matrices which have real determinant, and note that the
corresponding vector transformations~(\ref{2action}) are
determinant-preserving precisely when $\det(M)=\pm1$.

We are finally ready to define the octonionic transformation groups by
generalizing~(\ref{tsldef}), noting that the composition of linear
transformations is associative even when the underlying matrices are not
(since the order of operation is fixed).  However, in order to generate the
entire group, (compatible) transformations must be \textit{nested}; the action
of a composition of transformations can not in general be represented by a
single transformation.  We therefore generalize~(\ref{tsldef}) by defining
\begin{equation}
\tsl(2;\OO) := \Bigl\langle\left\{
	T_M : M\in\mm_2(\OO), \det(M)=\pm1 \right\}\Bigr\rangle
\label{tsldef2}
\end{equation}
where the angled brackets denote the span of the listed elements under
composition, and it is of course then straightforward to verify that
$\tsl(2;\OO)$ is a group under composition.  A similar definition can be given
for the spinor transformations, namely
\begin{equation}
\ssl(2;\OO)
  := \Bigl\langle\left\{
	S_M : M\in \mm_2(\OO), \det(M)=\pm1 \right\}\Bigr\rangle
\label{spsldef2}
\end{equation}

Since each transformation in $\tsl(2;\OO)$ preserves the determinant of
elements of $\JJ_2(\OO)$, it is clearly (isomorphic to) a subgroup of
$SO(9,1)$.  Manogue and Schray~\cite{Lorentz} showed, in slightly different
language, that in fact
\begin{equation}
\tsl(2;\OO) \cong SO(9,1)
\end{equation}
by giving an explicit set of basis elements which correspond to the standard
rotations and boosts in $SO(9,1)$.  Furthermore, it is easy to see that the
map
\begin{eqnarray}
\ssl(2;\OO) &\longrightarrow& \tsl(2;\OO) \\\nonumber
S_M &\longmapsto& T_M
\end{eqnarray}
is a two-to-one homomorphism with kernel $\{S_{\pm I}\}$, which establishes
the double covers
\begin{eqnarray}
	SL(2;\OO) &\longrightarrow& SO(9,1) \\
	SU(2;\OO) &\longrightarrow& SO(9)
\end{eqnarray}
(where $SU(2;\OO)$ is defined as for $SU(2;\HH)$ by restricting to
trace-preserving transformations), which are known results usually stated at
the Lie algebra level.

Despite the separate definitions presented above for $SL(2;\CC)$, $SL(2;\HH)$,
and $SL(2;\OO)$, a uniform definition can be given for any division algebra
$\KK=\RR,\CC,\HH,\OO$, modeled on the definition over $\OO$.  The basis used
by Manogue and Schray~\cite{Lorentz} consists of only two types of
transformations: single transformations corresponding to matrices of
determinant $+1$, and compositions of two transformations, each corresponding
to matrices of determinant $-1$; in this sense, each basis transformation can
be thought of as being ``of determinant $+1$''.  If we now define
\begin{eqnarray}
\ssla(2;\KK) &:=& \left\{S_M : M\in \mm_2(\KK),\ \det M=+1\right\} \nonumber\\
\sslb(2;\KK)
  &:=& \left\{S_P\circ S_Q : P,Q\in \mm_2(\KK),\ \det P=-1=\det Q \right\} \\
\ssl(2;\KK) &:=& \Bigl\langle \ssla(2;\KK) \cup \sslb(2;\KK) \Bigr\rangle
  \nonumber
\end{eqnarray}
where $\mm_2(\KK)$ denotes the set of \textit{complex} $2\times2$ matrices
over $\KK$, we recover the above definitions when $\KK=\HH,\OO$, while
retaining agreement with the standard definitions when $\KK=\RR,\CC$ (under
the usual identification of matrices with linear transformations).  A similar
definition can be made for $SU(2;\KK)$ by restricting to trace-preserving
transformations.

We can extend this treatment to the higher rank groups: There is a natural
action of $SL(n;\CC)$ as determinant-preserving linear transformations of
\nxn\ Hermitian (complex) matrices, with the unitary matrices $SU(n;\CC)$
additionally preserving the trace of \nxn\ Hermitian (complex) matrices, since
\begin{equation}
\tr(MXM^\dagger) = \tr(M^\dagger MX)
\label{traceid}
\end{equation}
and $M^\dagger M=I$ for $M\in SU(n;\CC)$, and these groups could be defined as
(the covering groups of) those groups of transformations.  Analogous results
hold for $SL(n;\HH)$ and $SU(n;\HH)$ (and of course also for $SL(n;\RR)$ and
$SU(n;\RR)$).

When extending these results to octonionic Hermitian matrices, we consider
only the $2\times2$ case discussed above and the $3\times3$ case, constituting
the \textit{exceptional Jordan algebra} $\JJ_3(\OO)$, also known as the
\textit{Albert algebra}.  In both cases, the determinant is well defined (see
below).  The group preserving the determinant in the $3\times3$ case is known
to be (a particular noncompact real form of) $E_6$; we can interpret this as
\begin{equation}
E_6 := \tsl(3;\OO) \cong SL(3;\OO)
\label{E6}
\end{equation}
Furthermore, the identity~(\ref{traceid}) from the complex case still holds
for $\XX\in\JJ_3(\OO)$ (and suitable $\MM\in E_6$, as discussed below) in the
form
\begin{equation}
\tr(\MM\XX\MM^\dagger) = \Re\Big(\tr(\MM^\dagger\MM\XX)\Big)
\end{equation}
where the right-hand side reduces to $\tr(\XX)$ if $\MM^\dagger\MM=\II$.  The
group which preserves the trace of matrices in $\JJ_3(\OO)$ is just (the
compact real form of)~$F_4$~\cite{Harvey}, which we can interpret as
\begin{equation}
F_4 \cong SU(3;\OO)
\label{F4}
\end{equation}
(There is no double-cover involved in~(\ref{E6}) and~(\ref{F4}), since these
real forms are simply-connected.)
At the Lie algebra level, this has been explained by Sudbery~\cite{Sudbery}
and at the group level this has been discussed by Ramond~\cite{Ramond} and
Freudenthal~\cite{Freudenthal}.

The remainder of this paper uses the above results from the $2\times2$ case to
provide an explicit construction of both $F_4$ and $E_6$ at the group level,
and discusses their properties.

\section{\boldmath Generators of $E_6$}

We consider octonionic $3\times3$ matrices $\MM$ acting on octonionic
Hermitian $3\times3$ matrices $\XX$, henceforth called \textit{Jordan
matrices}, in analogy with~(\ref{2action}), that is
\begin{equation}
T_\MM(\XX) = \MM\XX\MM^\dagger
\label{3action}
\end{equation}
For this to be well-defined, $\MM\XX\MM^\dagger$ must be Hermitian and hence
independent of the order of multiplication.  Just as was noted by Manogue and
Schray~\cite{Lorentz} in the $2\times2$ case, the necessary and sufficient
conditions for this are either that $\MM$ be complex or that the columns of
the imaginary part of $\MM$ be (real) multiples of each other.  As with
$SL(2;\OO)$, we will restrict ourselves to the case where $\MM$ is complex;
this suffices to generate all of $E_6$.

\subsection{Jordan Matrices}

The Jordan matrices form the exceptional Jordan algebra $\JJ_3(\OO)$ under the
commutative (but not associative) \textit{Jordan product} (see
e.g.~\cite{Harvey, Schafer})
\begin{equation}
\XX \circ \YY = {1\over2} (\XX\YY + \YY\XX)
\end{equation}
The \textit{Freudenthal product} of two Jordan matrices is given by
\begin{equation}
\XX * \YY = \XX \circ \YY - {1\over2} (\XX \, \tr(\YY) + \YY \, \tr(\XX))
	-{1\over2} \Big(\tr(\XX \circ \YY) - \tr(\XX) \, \tr(\YY)\Big)
\end{equation}
where the identity matrix is implicit in the last term.  The \textit{triple
product} of 3 Jordan matrices is defined by
\begin{equation}
[\XX,\YY,\ZZ] = (\XX * \YY) \circ \ZZ
\end{equation}
Finally, the \textit{determinant} of a Jordan matrix is defined by
\begin{equation}
\det \XX = \frac13 \> \tr [\XX,\XX,\XX]
\end{equation}
Remarkably, Jordan matrices satisfy the usual characteristic equation
\begin{equation}
\XX^3 - (\tr \XX) \, \XX^2 + \sigma(\XX) \, \XX - (\det \XX) \, \II = 0
\end{equation}
where we must be careful to define
\begin{equation}
\XX^3 := \XX^2 \circ \XX \equiv \XX \circ \XX^2
\end{equation}
and where the coefficient $\sigma(\XX)$ is given by
\begin{equation}
\sigma(\XX) := \tr(\XX*\XX) = {1\over2} \Bigl( (\tr \XX)^2 - \tr (\XX^2) \Bigr)
\end{equation}

\subsection{\boldmath $SO(9,1)$}

Consider first matrices of the form
\begin{equation}
\MM = \pmatrix{M& 0\cr \noalign{\smallskip} 0& 1\cr}
\label{Midentify}
\end{equation}
where $M\in SL(2;\OO)$ is one of the generators given by Manogue and
Schray~\cite{Lorentz}.  These generators include straightforward
generalizations of the standard representation of $SL(2;\CC)$ in terms of 3
rotations and 3 boosts, yielding 15 rotations and 9 boosts, together with
particular nested phase transformations (imaginary multiples of the identity
matrix), yielding the remaining 21 rotations corresponding to rotations of the
imaginary units ($SO(7))$.  Each such generator is complex and has real
determinant; for further details, and an explicit list of generators,
see~\cite{Lorentz}.  Since
\begin{equation}
\MM\pmatrix{X&\theta\cr \noalign{\smallskip} \theta^\dagger& n\cr}\MM^\dagger
  = \pmatrix{MXM^\dagger& M\theta\cr
	\noalign{\smallskip}
	(M\theta)^\dagger& n\cr}
  = \pmatrix{T_M(X)& S_M(\theta)\cr
	\noalign{\smallskip}
	\bigl(S_M(\theta)\bigr)^\dagger& n\cr}
\label{3rep}
\end{equation}
and using the fact that
\begin{equation}
\det\pmatrix{X&\theta\cr \noalign{\smallskip} \theta^\dagger& n\cr}
  = (\det X) n + 2 X \cdot \theta\theta^\dagger
\end{equation}
where
\begin{equation}
X\cdot Y = \frac12 \Bigl( \tr(X \circ Y) - \tr(X)\,\tr(Y) \Bigr)
\end{equation}
is the Lorentzian inner product in $9+1$ dimensions, it is straightforward to
verify that $T_\MM$ preserves the determinant of a Jordan matrix $\XX$, and is
hence in $E_6$, if $\MM$ is of the form~(\ref{Midentify}).  This shows that
\begin{equation}
  SL(2;\OO) \subset E_6
\end{equation}
as expected, where, as already noted, $SL(2;\OO)$ is the double cover of
$SO(9,1)$.  Since $SL(2;\OO)$ acts as $SO(9,1)$ on $X$ in~(\ref{3rep}), we
will somewhat loosely describe $SO(9,1)$ itself (and its subgroups) as being
``in'' $E_6$.

This construction in fact yields three obvious copies of $SO(9,1)$ contained
in $E_6$, corresponding to the three natural ways of embedding a $2\times2$
matrix inside a $3\times3$ matrix.  These three copies are related by the
\textit{cyclic permutation matrix}
\begin{equation}
\TT = \pmatrix{0& 1& 0\cr 0& 0& 1\cr 1& 0& 0\cr}
\end{equation}
which satisfies
\begin{equation}
\TT^{-1} = \TT^2 = \TT^\dagger
\end{equation}
and which is clearly in $E_6$.

Conversely, all elements of $E_6$ can be built up out of these (three sets of)
$SO(9,1)$ transformations.

\subsection{\boldmath $SO(8)$, Triality, and $F_4$}

Since each copy of $SO(9,1)$ is 45-dimensional, but the dimension of $E_6$ is
only 78, it is clear that our description so far must contain some redundancy.
We note first of all that~(\ref{3rep}) contains not only the vector
representation~(\ref{2action}) of $SO(9,1)$, but also the dual spinor
representations
\begin{eqnarray}
\theta &\longmapsto& M\theta \\
\theta^\dagger &\longmapsto& \theta^\dagger M^\dagger \nonumber
\end{eqnarray}
and therefore combines $2\times2$ vector and spinor representations into a
single $3\times3$ representation.  \textit{Triality} says that, for $SO(8)$,
these three representations are isomorphic.
\footnote{The term \textit{triality} appears to have first been used by
Cartan~\cite{Cartan}, who used it to describe the symmetries of the Dynkin
diagram of $SO(8)$.  An infinitesimal \textit{principle of triality} in the
language of derivations is proved in~\cite{Schafer}, which  credits
Jacobson~\cite{JacobsonII} with the analogous theorem for Lie groups.
Baez~\cite{Baez} describes the four \textit{normed trialities} as trilinear
maps on representations of (particular) Lie algebras, and discusses their
relationship to the four normed division algebras and their automorphisms.  A
similar treatment can be found in Conway and Smith~\cite{Conway}.}

To see explicitly why triality holds, we begin with the description of
$SO(8)\subset SO(9,1)$ from~\cite{Lorentz}.  Since $SO(8)$ transformations of
the form~(\ref{2action}) leave the diagonal of $X$ invariant, these
transformations correspond to the ``transverse'' degrees of freedom in
$SO(9,1)$.  One might therefore expect $SO(8)$ transformations to take the
form $\pmatrix{q& 0\cr 0& \bar{r}\cr}$ with $|q|=|r|=1$.  However, the
essential insight of~\cite{Lorentz} was to require that all $SO(9,1)$
transformations be \textit{compatible}, that is, that they (be generated by
matrices which) satisfy~(\ref{comp}); we will see the importance of this
requirement below.  This condition restricts the allowed form of $SO(8)$
transformations to those which can be constructed from ($2\times2$) diagonal
matrices which are either imaginary multiples of the identity matrix, or of
the form
\begin{equation}
M = \pmatrix{q& 0\cr 0& \bar{q}\cr}
\label{so8}
\end{equation}
where $|q|=1$, so that $q=e^{s\theta}$ for some imaginary unit $s\in\OO$ with
$s^2=-1$.  As discussed in~\cite{Lorentz}, the matrix~(\ref{so8}) induces a
rotation in the ($1,s$)-plane through an angle $2\theta$.  Furthermore,
$SO(7)$ transformations, namely those leaving invariant the identify element
$1$, can be constructed by suitably nesting an even number of purely imaginary
matrices of the form~(\ref{so8}), that is, matrices of this form for which
$\theta={\pi\over2}$.  This allows us to generate all of $SO(8)$ using
matrices which have determinant~$1$.  Alternatively, $SO(7)$ transformations
can also be obtained by nesting imaginary multiples of the identity matrix
(which have determinant~$-1$), since this involves an even number of sign
changes when compared with the above description.

Inserting~(\ref{so8}) into~(\ref{Midentify}), the resulting $E_6$
transformation $\MM$ leaves the diagonal of a Jordan matrix $\XX\in\JJ_3(\OO)$
invariant.  Explicitly, writing
\begin{equation}
\XX = \pmatrix{p& \bar{a}& c\cr a& m& \bar{b}\cr \bar{c}& b& n\cr}
\end{equation}
we see that the action~(\ref{3action}) leaves $p$, $m$, $n$ invariant, and
acts on the octonions $a$, $b$, $c$ via
\begin{eqnarray}
a &\longmapsto& \bar{q}a\bar{q} \nonumber\\
b &\longmapsto& bq \label{so8action}\\
c &\longmapsto& qc \nonumber
\end{eqnarray}
These three transformations are precisely the standard description of the
(vector and two spinor) representations of $SO(8)$ in terms of symmetric,
left, and right multiplication by unit octonions.  The
actions~(\ref{so8action}) provide an implicit mapping between these three
representations (obtained by using the same $q$ in each case),
\footnote{To verify this assertion, one must first argue that the implicit map
between representations in~(\ref{so8action}) is well-defined.  At issue is the
uniqueness of representations such as $a \longmapsto q_1(...(q_m a)...)$ for a
particular $SO(8)$ transformation, and specifically whether this notion of
uniqueness is the same for the three representations.  One way to show this is
to explicitly construct the maps between the representations.}
which is clearly both a (local) diffeomorphism and a 1-to-1 map between the
two spinor representations, and a 2-to-1 map between either spinor
representation and the vector representation.  For us, \textit{triality} is
this explicit relationship between the three representations.

In our language, this means that if $M$ in~(\ref{Midentify}) is an $SO(8)$
transformation, then not only does $\MM$ generate an $SO(8)$ transformation on
$\XX$ via~(\ref{3rep}), but so do $\TT\MM\TT^2$ and $\TT^2\MM\TT$, since each
of these latter two transformations differs from $\MM$ merely in which
representation of $SO(8)$ acts on each of $a$, $b$, $c$.  Even though these
individual transformations are different, the collection of all of them is the
same in each case.  Note that this identification of the three copies of
$SO(8)$ is only possible because the original $SO(9,1)$ transformation was
assumed to be compatible.  Thus, (the double cover of) $SO(8)$ is precisely
the subgroup of $E_6$ which leaves the diagonal of every Jordan matrix $\XX$
invariant.

The dimension of the single resulting copy of $SO(8)$ is 28.  Adding in the
$3\times8=24$ additional rotations in (the three copies of) $SO(9)$ yields
$52$, the correct dimension for $F_4$.  Including the $3\times9-1=26$
independent boosts gives the full $78$ generators of $E_6$.  Thus, triality
fully explains the redundancy in our original 135 generators.

At the Lie algebra level, the dimension of $E_6$ can be determined by first
noting that the diagonal elements are not independent.  It turns out that an
independent set can be taken to be the 64 independent tracefree matrices,
together with the 14 generators of $G_2$ (see below), for a total of 78
generators.  Of these, 26 (24 non-diagonal + 2 diagonal) are Hermitian and
hence boosts (the third diagonal Hermitian generator is not independent); the
remaining 52 generators yield $F_4$.
In fact, the (complex) generators of $SO(9,1)$ as given in~\cite{Lorentz} all
satisfy either $MM^\dagger=I$ (rotations) or $M=M^\dagger$ (boosts).  But
$F_4$ is generated precisely by the unitary elements of $E_6$, and hence is
generated by the (3 sets of) rotations in $SO(9)$.

\subsection{\boldmath $G_2$}

As discussed in~\cite{Lorentz}, the automorphism group of the octonions,
$G_2$, can be constructed by suitably combining rotations of the octonionic
units, thus providing explicit verification that $G_2$ is a subgroup of
$SO(7)$.  In particular, a copy of $G_2$ sits naturally inside each $SO(9,1)$,
generated by 14 (nested) imaginary multiples of the identity matrix (the
``additional transverse rotations'' of~\cite{Lorentz}, also denoted
``flips'').
\footnote{Even though these three copies of $SO(7)$ all live in the single
copy of $SO(8)$ described above, they are not the same.}
We thus appear to have three copies of $G_2$ sitting inside $E_6$, one for
each copy of $SO(9,1)$.

As further shown in~\cite{Lorentz}, the automorphisms of $\OO$ can be
generated by octonions of the form
\begin{equation}
e^{\hat{q}\theta} = \cos(\theta) + \hat{q} \sin(\theta)
\end{equation}
with $\hat{q}$ a pure imaginary, unit octonion, but where $\theta$ must be
restricted to be a multiple of $\pi/3$, corresponding to the sixth roots of
unity.  But, as can be verified by direct computation, multiplying the
identity matrix by such an automorphism leads to an element of $E_6$, thus
giving us yet another apparent copy of $G_2$ in $E_6$.

Remarkably, due to triality, all four of these subgroups are the same.

To see this, consider the rotations by $\pi\over2$ (``flips'') used
in~\cite{Lorentz} to generate the transverse rotations.  Using the
identification (\ref{Midentify}), such a transformation takes the form
\begin{equation}
\QQ_{\hat{q}} = \pmatrix{\hat{q}I& 0\cr \noalign{\smallskip} 0& 1\cr}
\end{equation}
and we have
\begin{equation}
\QQ_{\hat{q}}	
	\pmatrix{X&\theta\cr \noalign{\smallskip} \theta^\dagger& n\cr}
	\QQ_{\hat{q}}^\dagger
  = \pmatrix{-\hat{q}X\hat{q}&\hat{q}\theta\cr
	\noalign{\smallskip} -\theta^\dagger\hat{q}& n\cr}
\end{equation}
Under this transformation, $X$, $\theta$, and $\theta^\dagger$ undergo
separate $SO(7)$ transformations, related by triality.  We emphasize that, in
general, the off-diagonal elements of $\XX$ undergo \textit{different} $SO(7)$
transformations.  (The diagonal elements are of course fixed by any such
transformation.)

Acting on a single octonion, nested sequences of these $SO(7)$ transformations
can be used to generate $G_2$.  For instance, conjugating successively with
\begin{equation}
\hat{q} = i, i\cos\theta+i\ell\sin\theta, j, j\cos\theta-j\ell\sin\theta
\end{equation}
yields a $G_2$ transformation which leaves the quaternionic subalgebra
generated by $k$ and $\ell$ fixed.  What happens when this sequence of
$\hat{q}$'s is applied as $E_6$ transformations, that is, in the form
$\QQ_{\hat{q}}$?  Remarkably, direct computation shows in this case that the
elements of $\XX$ all undergo the same $G_2$ transformation.  Since all $G_2$
transformations can be generated by such transformations, triality is, in this
sense, the identity map on $G_2$!  The $G_2$ transformation obtained by
suitably nesting $\QQ_{\hat q}$'s is therefore the same as the $G_2$
transformation obtained by replacing $\QQ_{\hat q}$ by $\hat{q}\II$ at each
step.  This shows that the three $G_2$ subgroups contained in the three copies
of $SO(9,1)$ are all identical to the ``diagonal'' $G_2$ subgroup, as claimed
above.

An explicit example of triality-related automorphisms is given by
\footnote{Further examples can be found in~\cite{Mobius}.}
\begin{equation}
k(j(iq))) = k(j(iq\ibar)\jbar)\kbar = (((q\ibar)\jbar)\kbar)
\end{equation}
with $q\in\OO$, which realizes ``$\ell$-conjugation'' as a linear map.

\section{\Spinor s}

We have argued elsewhere~\cite{Dim,Spin} that the ordinary momentum-space
(massless and massive) Dirac equation in $3+1$ dimensions can be obtained via
dimensional reduction from the Weyl (massless Dirac) equation in $9+1$
dimensions.  The dimensional reduction is accomplished by the simple expedient
of choosing a preferred complex subalgebra of the octonions, thus reducing
$SL(2;\OO)$ to $SL(2;\CC)$, and hence the Lorentz group in 10 spacetime
dimensions to that in 4 dimensions.

The massless Dirac equation in 10 spacetime dimensions can be written in
momentum space as the eigenvalue problem
\begin{equation}
\label{Dirac}
\tilde{P}\psi = 0
\end{equation}
where $P$ is a $2\times2$ octonionic Hermitian matrix corresponding to the
10-dimensional momentum vector, $\psi\in\OO^2$ is a 2-component octonionic
column, corresponding to a Majorana-Weyl spinor, and where tilde denotes trace
reversal, that is
\begin{equation}
\tilde{P} = P - \tr(P)\,I
\end{equation}
The general solution of~(\ref{Dirac}) is
\begin{eqnarray}
P &=& \pm\theta\theta^\dagger \\
\psi &=& \theta\xi
\end{eqnarray}
where $\theta\in\OO^2$ is a 2-component octonionic vector whose components lie
in the same complex subalgebra of $\OO$ as do those of $P$, and where
$\xi\in\OO$ is arbitrary.  (Such a $\theta$ must exist since $\det(P)=0$.)

In~\cite{Spin}, we further showed how to translate the standard treatment of
the Dirac equation in terms of gamma matrices into octonionic language,
pointing out that a 2-component \textit{quaternionic} formalism is of course
isomorphic to the traditional 4-component complex formalism.  Remarkably, the
above solutions to the \textit{octonionic} Dirac equation must be
quaternionic, as they only involve 2 independent octonionic directions.  This
allows solutions of the octonionic Dirac equation to be interpreted as
standard fermions --- and one can fit precisely 3 ``families'' of such
quaternionic solutions into the octonions, which we interpret as generations.
For further details, see~\cite{Spin}, or the more recent treatment
in~\cite{York}.

As outlined in~\cite{Other}, it is natural to introduce a 3-component
formalism; this approach was first suggested to us by Fairlie and
Corrigan~\cite{FC}, and later used by Schray~\cite{Schray,JorgThesis} for the
superparticle.  Defining
\begin{equation}
\Psi = \pmatrix{\theta\cr \bar\xi}
\end{equation}
we have first of all that
\begin{equation}
\PP := \Psi \Psi^\dagger = \pmatrix{P& \psi\cr \psi^\dagger& |\xi|^2\cr}
\label{1sq}
\end{equation}
so that $\Psi$ combines the bosonic and fermionic degrees of freedom.  Lorentz
transformations on both the vector $P$ and the spinor $\psi$ now take the
elegant form~(\ref{3rep}), which we used to view $SO(9,1)$ as a subgroup of
$E_6$; the rotation subgroup $SO(9)$ lies in $F_4$.  We refer to $\Psi$ as a
\textit{\spinor}.

Direct computation shows that the Dirac equation~(\ref{Dirac}) is equivalent
to the equation
\begin{equation}
\PP*\PP=0
\label{DiracF}
\end{equation}
whose solutions are precisely \textit{quaternionic} matrices of the
form~(\ref{1sq}), that is, (the components of) $\theta$ and $\xi$ must lie in
a quaternionic subalgebra of $\OO$; $\Psi$ is a quaternionic \spinor.  But the
\textit{Cayley plane} $\OO{\mathbb P}^2$ consists of those elements
$\PP\in\JJ_3(\OO)$ which satisfy
\begin{equation}
\PP \circ \PP = \PP; \qquad \tr\,\PP = 1
\end{equation}
and this turns out to be equivalent to requiring $\PP$ to be a (normalized)
solution of~(\ref{DiracF}).  Thus, in this interpretation, the Cayley plane
consists precisely of normalized, quaternionic \spinor s, and these are
precisely the (normalized) solutions of the Dirac equation.

Furthermore, any Jordan matrix can be decomposed in the form~\cite{Paige}
\begin{equation} 
\AA = \sum_{i=1}^3 \lambda_i \PP_i
\label{Decomp}
\end{equation}
with
\begin{equation}
\PP_i\circ \PP_j = 0 \qquad(i\ne j)
\end{equation}
As discussed in~\cite{Other}, we refer to the decomposition~(\ref{Decomp}) as
a \textit{$p$-square decomposition} of $\AA$, with $p$ denoting the number of
nonzero eigenvalues $\lambda_i$.  The $\PP_i$ are eigenmatrices of $\AA$, with
eigenvalue $\lambda_i$, that is
\begin{equation}
\AA \circ \PP_i = \lambda_i \PP_i
\end{equation}
and $p$ denotes the number of nonzero eigenvalues, and hence the number of
nonzero primitive idempotents in the decomposition.  As shown in~\cite{Other},
if $\det(\AA)\ne0$, then $\AA$ is a 3-square, while if
$\det(\AA)=0\ne\sigma(\AA)$, then $\AA$ is a 2-square.  Finally, if
$\det(\AA)=0=\sigma(\AA)$, then $\AA$ is a 1-square (unless also $\tr(\AA)=0$,
in which case $\AA\equiv0$).  It is intriguing that, since $E_6$ preserves
both the determinant and the condition $\sigma(\AA)=0$, $E_6$ therefore
preserves the class of $p$-squares for each $p$.  But solutions of the Dirac
equation~(\ref{DiracF}) are 1-squares!  Thus, the Dirac equation in 10
dimensions admits $E_6$ as a symmetry group.

The particle interpretation described in~\cite{Dim,Spin} suggests regarding
1-squares as representing three generations of leptons.  If 1-squares
correspond to leptons, could it be that 2-squares are mesons and 3-squares are
baryons?

\section{The Structure of \boldmath$E_6$}

We have shown that the massless 10-dimensional Dirac equation, originally
posed as an eigenvalue problem for $2\times2$ octonionic Hermitian matrices,
is in fact equivalent to the defining condition for the Cayley plane.  This
suggests that the natural arena for the Dirac equation is a 3-component
formalism involving Cayley spinors, which explicitly incorporates both
bosonic and fermionic degrees of freedom, suggesting a natural supersymmetry.
Furthermore, the symmetry group of the Dirac equation has been shown to be
$E_6$, suggesting that $E_6$ (or possibly one of its larger cousins, $E_7$ or
$E_8$) is the natural symmetry group of fundamental particles.

\begin{figure}
\epsfysize=5truein
\epsffile[84 353 478 685]{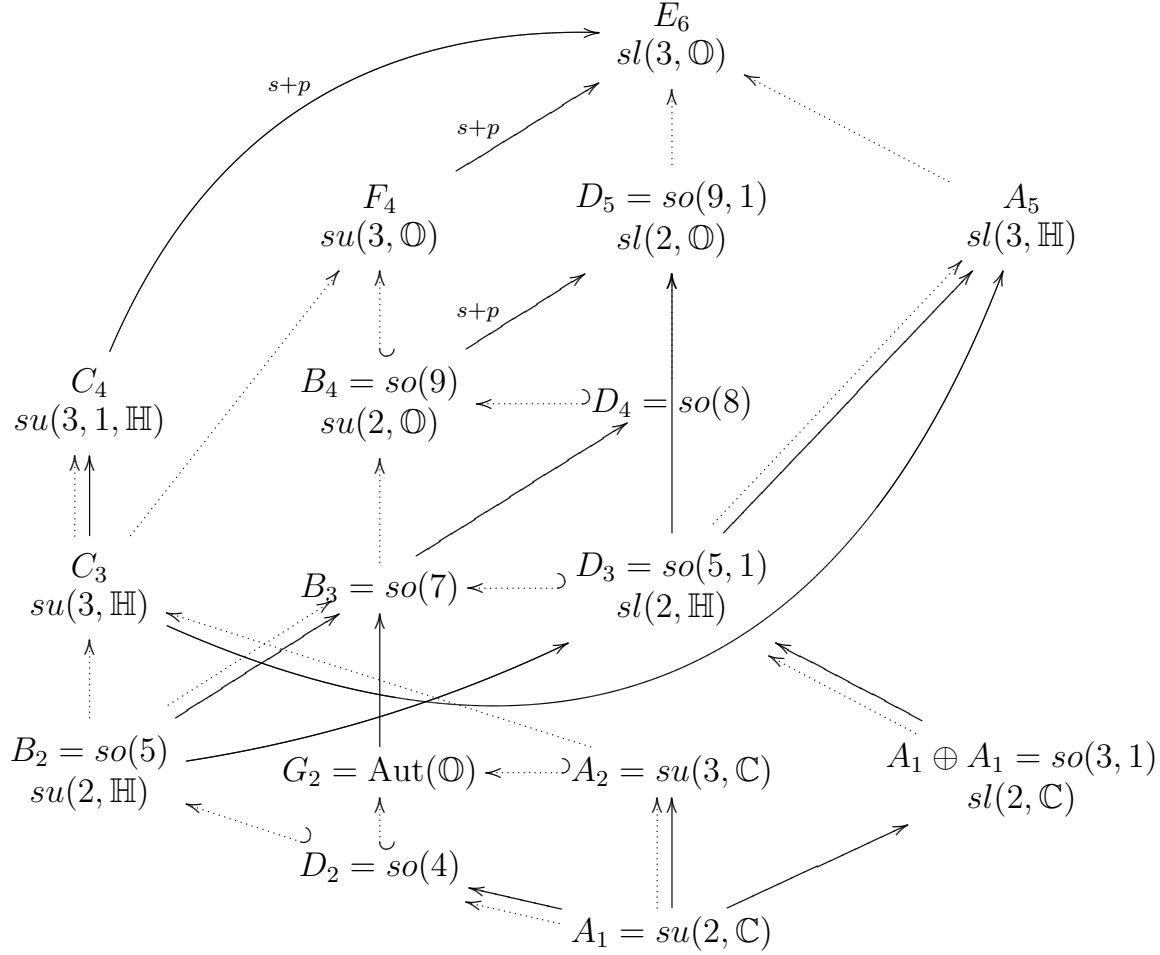}
\caption{A map of $E_6$ (taken from~\cite{AaronThesis}).}
\label{E6map}
\end{figure}

Understanding the structure of (this particular real representation of) $E_6$
may therefore be of great importance to an ultimate understanding of
fundamental particles.  In this regard, we call the reader's attention to the
recent work of Aaron Wangberg~\cite{AaronThesis,JOMA}, which describes the
real representations of $E_6$ and its physically important subgroups.  A
``map'' of $E_6$, excerpted from~\cite{AaronThesis}, appears in
Figure~\ref{E6map}.

\section*{Acknowledgments}

This work was supported in part by a grant from the Foundational Questions
Institute (FQXi).  We thank the Department of Physics at Utah State University
for support and hospitality during the preparation of this manuscript.

\end{document}